\documentclass[11pt]{article}
\usepackage{amssymb}
\usepackage{latexsym}

\newtheorem{thm}{Theorem}[section]
\newtheorem{lem}[thm]{Lemma}
\newtheorem{cor}[thm]{Corollary}
\newtheorem{pro}[thm]{Proposition}

\newtheorem{defi}[thm]{Definition}

\newcommand{\gm }{\Gamma }
\newcommand{\lon }{\longrightarrow }
\newcommand{\be }{\begin{eqnarray*}}
\newcommand{\ee }{\end{eqnarray*}}

\newcommand{\pf}{\noindent{\bf Proof.}\ }
\newcommand{\qed}{\begin{flushright} $\Box$\ \ \ \ \ \
                    \end{flushright}}
\newcommand{\complex}{{\Bbb C}}
\newcommand{\reals}{{\Bbb R}}

\newcommand{\frakg}{{\frak g}}
\newcommand{\frakh}{{\frak h}}
\newcommand{\frakl}{{\frak l}}
\newcommand{\frakm}{{\frak m}}

\newcommand{\hstar}{*_{\hbar}}

\newcommand{\half}{\frac{1}{2}}

\newcommand{\cald}{{\cal D}}

\def\description label#1{\hfil\bf[#1]\hfil}
\newcommand{\parr}[1]{\frac{\partial  #1}{\partial \lambda_{i}}}

\newcommand{\ot}{\mbox{$\otimes$}}

\newcommand{\alt}{\mbox{Alt}}

\newcommand{\td}{{\Delta}}
\newcommand{\ttd}{\widetilde{\Delta}}
\newcommand{\vh}[1]{\Vec{h_{#1}}}
\newcommand{\vhh}[1]{\Vec{h_{#1}'}}
\newcommand{\sstar}{*}
\newcommand{\prefer}{compatible }
\newcommand{\flb}{{[\![}}
\newcommand{\frb}{{]\!]}}   

\def\sdp{\mathbin{\hbox{$\mapstochar\kern-.3333em\times$}}}
\def\pds{\mathbin{\hbox{$\times\kern-.55em\mapstochar\,$}}}

\newcommand{\wed}{\mathbin{\lower1.5pt\hbox{$\scriptstyle{\wedge}$}}}

\let\Tilde=\widetilde

\let\Vec=\overrightarrow

\def\chigh{{\raise1.5pt\hbox{$\chi$}}}
\let\phi=\varphi
\def\til0{\Tilde{0}}

\def\dminus{\raise2pt\hbox{\vrule height1pt width 2ex}\hskip3pt}

\def\pback#1{\mathbin{{{\lower1.2ex\hbox{$\times$}}\atop #1}}}

\def\vlra{\hbox{$\,-\!\!\!-\!\!\!-\!\!\!-\!\!\!-\!\!\!
-\!\!\!-\!\!\!-\!\!\!-\!\!\!-\!\!\!\longrightarrow\,$}}

\def\gpd{\,\lower1pt\hbox{$\longrightarrow$}\hskip-.24in\raise2pt
               \hbox{$\longrightarrow$}\,}

\def\lgpd{\,\lower1pt\hbox{$\vlra$}\hskip-1.02in\raise2pt\hbox{$\vlra$}\,}

\def\llgpd{\,\lower1pt\hbox{$\vvlra$}\hskip-1.3in\raise2pt\hbox{$\vvlra$}\,
}

\begin{document}

\title{{\bf Quantum dynamical Yang-Baxter equation over a nonabelian base}}
\author{ PING XU \thanks{ Research partially supported by   NSF
          grant DMS00-72171.}\\
   Department of Mathematics\\The  Pennsylvania State University \\
University Park, PA 16802, USA\\
          {\sf email: ping@math.psu.edu }}

\date{}

\maketitle

\begin{abstract} 
In this paper we consider dynamical r-matrices over a
nonabelian base. There are two main results. First,  
corresponding to a fat reductive decomposition of a Lie algebra $\frakg
=\frakh \oplus \frakm$,  we   construct geometrically
a non-degenerate triangular dynamical
r-matrix using  symplectic
fibrations.  Second,  we  prove that a triangular
dynamical r-matrix $r: \frakh^* \lon \wedge^2 \frakg$
corresponds to a Poisson manifold $\frakh^* \times G$.
A special type of quantizations of this Poisson manifold, called  \prefer
star products in this  paper, yields  
a generalized version of the quantum dynamical Yang-Baxter    equation
 (or Gervais-Neveu-Felder equation).  As a result, the
quantization problem of a general dynamical r-matrix
is proposed.
\end{abstract}

\section{Introduction}

Recently, there  has been  growing interest in the so called quantum dynamical
Yang-Baxter equation:
\begin{equation}
\label{eq:dybe0}
R_{12}(\lambda )  R_{13}(\lambda +\hbar h^{(2)} ) 
  R_{23}(\lambda )
=R_{23}(\lambda +\hbar h^{(1)} )  
R_{13}(\lambda )   R_{12}(\lambda +\hbar h^{(3)} ).
\end{equation}
This equation  arises
naturally from  various contexts in mathematical physics.
It first appeared in the work of Gervais-Neveu in their
study of quantum Liouville theory \cite{GN}.  Recently it reappeared
in Felder's work on the  quantum Knizhnik-Zamolodchikov-Bernard
 equation \cite{Felder}. It also  has  been found
to be connected  with the quantum Caloger-Moser systems
 \cite{ABE}.
As the  quantum Yang-Baxter equation is connected
with quantum groups, the  quantum dynamical
Yang-Baxter equation is known to be connected with elliptic
quantum groups \cite{Felder},
as well as with Hopf algebroids or quantum groupoids \cite{EV2, Xu2, Xu3}.

 The classical counterpart
of the quantum dynamical Yang-Baxter equation  was first
considered by Felder \cite{Felder}, and then studied
by Etingof and Varchenko \cite{EV1}. This is
 the so called classical dynamical
Yang-Baxter equation, and a solution to such an
equation (plus some other reasonable conditions) is
called a classical dynamical r-matrix. More precisely,
given a Lie algebra $\frakg$ over $\reals$ (or over $\complex$)
 with an Abelian Lie subalgebra $\frakh$, a classical dynamical r-matrix
is a smooth (or meromorphic)  function $r: \frakh^* \lon
\frakg \ot \frakg $
satisfying the following conditions:
\begin{enumerate}
\item (zero weight condition) $[h\ot 1 +1 \ot h, r(\lambda )]=0, \  \ 
\forall h\in \frakh$;
\item (normal condition)  $r_{12}+r_{21}= \Omega$, where
$\Omega\in (S^{2}\frakg)^{\frakg}$ is a Casimir element;
\item (classical dynamical Yang-Baxter equation\footnote{Throughout
the paper, we follow the sign convention in \cite{ABE} for  the
definition  of a  classical dynamical $r$-matrix in order to be
consistent with the quantum dynamical Yang-Baxter
equation (\ref{eq:dybe0}). This  differs  in a sign from
the one used in \cite{EV1}.})
\begin{equation}
\label{cdybe}
 Alt(dr) \, - \,
([r_{12},r_{13}] + [r_{12}, r_{23}] + [r_{13}, r_{23}]) \, = \, 0,
\end{equation}
\end{enumerate}
where $\alt dr =\sum (h_{i}^{(1)} \frac{\partial r_{23}}{\partial \lambda^{i}}
-h_{i}^{(2)} \frac{\partial r_{13}}{\partial \lambda^{i}}+
h_{i}^{(3)} \frac{\partial r_{12}}{\partial \lambda^{i}} )$.

 A fundamental question is whether a  classical dynamical r-matrix
is always quantizable.  There has appeared  a lot of
work in this direction, for example, see \cite{Arnaudon, Jimbo, ESS}.
 In the  triangular 
case (i.e.,  $r$ is skew-symmetric: $r_{12}(\lambda ) +r_{21}(\lambda ) =0$),
 a general quantization  scheme was developed by the author using
the Fedosov method, which works 
 for a vast class of dynamical
r-matrices, called splittable triangular dynamical 
r-matrices \cite{Xu4}. 
Recently, Etingof and  Nikshych, using  the vertex-IRF 
transformation method, proved  the existence of
 quantizations  for the so called completely degenerate
 triangular dynamical r-matrices \cite{EN}.

Interestingly, although the quantum  dynamical
 Yang-Baxter equation in \cite{Felder} only makes
sense when the base Lie algebra $\frakh$ is Abelian,
its classical counterpart admits an immediate generalization
for any base Lie algebra  $\frakh$ which is not necessarily
Abelian. Indeed, all one needs  to do is to change
 the first  condition (i) to:

(i'). $r: \frakh^* \lon \frakg \ot \frakg $ is $H$-equivariant,
where $H$ acts on $\frakh^*$ by coadjoint action and
on $\frakg \ot \frakg $ by adjoint action.

There exist many examples of such classical dynamical r-matrices.
For instance, when $\frakg$ is a simple Lie algebra and $\frakh$ is a reductive
Lie subalgebra containing the Cartan subalgebra,
there is a classification due to  Etingof-Varchenko  \cite{EV1}.
In particular, when  $\frakh=\frakg$,
an explicit formula was discovered by Alekseev and Meinrenken
in their study of non-commutative Weil algebras \cite{AM}. Later,
this was generalized by  Etingof and Schiffermann \cite{ES}
to a more general context. Moreover, under some regularity condition,
 they showed that the moduli space
of dynamical r-matrices  essentially consists of a single point once the
initial value of the dynamical r-matrices  is fixed.
 A natural question
arises as to what should be the quantum counterpart
of these r-matrices. And more generally, is any
classical dynamical r-matrix (with nonabelian base)
quantizable?

A basic question is  what the   quantum dynamical
Yang-Baxter equation should look like when $\frakh$ is nonabelian.
 In this paper, as a toy model,  we   consider the special case of triangular 
dynamical r-matrices and their quantizations.
 As in the Abelian case, 
these r-matrices  naturally   correspond to  some invariant
Poisson structures on $\frakh^* \times G$.
It is standard that quantizations of Poisson structures
correspond to star products \cite{BEFFL}.
The special form of the Poisson bracket relation on $\frakh^* \times G$ 
suggests a specific form that their star products should
take. This leads to our definition of \prefer
star products. The compatibility condition (which, in this
case,   is just the associativity)
naturally leads to a  quantum dynamical Yang-Baxter equation: 
Equation (\ref{eq:dybe}). As we shall  see, this equation  indeed
resembles  the   usual quantum dynamical Yang-Baxter equation
(unsymmetrized version). The only difference is
that  the usual pointwise  multiplication on $C^{\infty}(\frakh^* )$
is replaced by  the PBW-star product, which is  indeed
the deformation quantization of the canonical
Lie-Poisson structure on $\frakh^* $.
Although  Equation (\ref{eq:dybe}) is derived  by considering 
triangular dynamical r-matrices, it makes perfect
sense  for non-triangular ones as well.  This
naturally leads to our definition
of quantization of dynamical  r-matrices over
an arbitrary base Lie subalgebra which is not
necessary Abelian. The 
problem is that such an  equation only makes sense for
$R: \frakh^* \lon U\frakg \ot U\frakg \flb\hbar \frb$. In
the Abelian case, it appears that one may  consider $R$ valued in a
deformed universal enveloping algebra  $U_{\hbar} \frakg$,
but in most cases $U_{\hbar} \frakg$ is isomorphic
to $U\frakg \flb\hbar \frb$ as an algebra. So Equation (\ref{eq:dybe}),
 in a certain sense,  is general enough to include all the interesting cases. 
However, the physical meaning of this equation  remains mysterious.

Another main  result  of the paper is to give  a
geometric construction of triangular dynamical
r-matrices. More precisely, we give  an explicit 
construction of a triangular dynamical   r-matrix
from a fat reductive decomposition of a Lie algebra 
 $\frakg=\frakh\oplus \frakm$ (see  Section 2 for the definition).
This includes those
examples of triangular dynamical r-matrices considered
in \cite{EV1}. Our  main purpose is to show 
that triangular dynamical r-matrices (with nonabelian base)
do rise naturally from symplectic geometry. 
This gives us  another  reason why it is important to consider
their  quantizations.
Discussion on this part occupies Section 2. 
Section 3 is devoted to the discussion  of 
\prefer star products, whose associativity
leads to a ``twisted-cocycle" condition.
In Section 4, we will derive the quantum
dynamical Yang-Baxter equation from this
twisted-cocycle condition. The last section contains
some concluding  remarks and open questions. 

Finally, we note that in this paper,  by a 
dynamical r-matrix, we always  mean a 
dynamical r-matrix over a general base
Lie subalgebra unless specified. Also 
 Lie algebras are normally  assumed to be
over $\reals$, although most   results can
be  easily modified  for complex Lie algebras.

{\bf Acknowledgments.} The author  would   like to thank
 Philip Boalch, Pavel Etingof,  Boris Tsygan and David Vogan 
for fruitful discussions and comments. He is  especially grateful
to Pavel Etingof for explaining the paper \cite{ES}, which
inspired his interest on this topic.
He  also wishes  to thank  Simone Gutt and Stefan Waldmann
 for  providing  him   some useful references on star products
of  cotangent symplectic manifolds.

\section{Classical dynamical $r$-matrices}

In  this section, we will
 give a geometric construction of
 triangular dynamical r-matrices.
As we shall see, 
these r-matrices do  arise naturally from
symplectic geometry. We will show some
interesting examples, which  include triangular dynamical r-matrices
for simple Lie algebras constructed by Etingof-Varchenko  \cite{EV1}.


Below let us recall the definition
 of a  classical triangular dynamical $r$-matrix.
Let $\frakg$ be  a Lie algebra over $\reals$ (or $\complex$)
and $\frakh \subset \frakg$ a Lie subalgebra.
A classical dynamical r-matrix $r : \frakh^* \lon \frakg\ot \frakg$
is said to be {\em triangular} if it is skew symmetric:
$r_{12}+r_{21}=0$. In other words, a 
 classical triangular dynamical $r$-matrix  is a smooth function (or
meromorphic function in the complex case)
   $r: \frakh^{*} \lon \wedge^{2} \frakg$ such that
\begin{enumerate}
\item $r: \frakh^* \lon \wedge^{2} \frakg $  is $H$-equivariant,
where $H$ acts on $ \frakh^*$ by coadjoint action and 
acts on $\wedge^{2} \frakg $ by adjoint action.
\item 
\begin{equation}
\label{eq:triangular}
\sum_{i} h_{i} \wedge \frac{\partial r}{ \partial\lambda^{i}} 
-\half [r, r]=0,
\end{equation}
\end{enumerate}
where the bracket $[\cdot , \cdot ]$ refers to the Schouten type
bracket: $\wedge^{k}\frakg \ot \wedge^{l}\frakg \lon \wedge^{k+l-1}\frakg $
induced from the Lie algebra bracket on $\frakg$,
 $\{h_{1}, \cdots , h_{l}\}$  is a basis of  $\frakh$,  and
$(\lambda^{1}, \cdots , \lambda^{l})$ its  induced
coordinate system  on $\frakh^*$. 

The following proposition gives an alternative description of a
classical triangular dynamical r-matrix. 

\begin{pro}
\label{pro:dy}
A smooth function  $r:\frakh^{*} \lon \wedge^2 \frakg  $
 is a triangular dynamical $r$-matrix iff
$$\pi = \pi_{\frakh^* } + \sum_{i}  \parr{}\wedge \vh{i} +\Vec{r(\lambda )}$$
 is a  Poisson tensor
 on $M=\frakh^* \times G$, where  $\pi_{\frakh^* }$ denotes
the standard Lie  (also known as Kirillov-Kostant)
Poisson tensor on the Lie algebra dual $\frakh^* $,
  $\vh{i}\in {\frak X} (M) $ is the 
left invariant vector field on $M$ generated by $h_{i} \in \frakh$,
 and  similarly $\Vec{r(\lambda )}\in \gm (\wedge^{2} TM)$ is the
 left invariant bivector field
on $M$ corresponding to  $r(\lambda )$.
\end{pro}
\pf  Set $$\pi_{1}=\pi_{\frakh^* } + \sum_{i}  \parr{}\wedge \vh{i}.$$
Then $\pi =\pi_{1}+\Vec{r(\lambda )}$. Note that, for any $(\lambda , x)$, 
$\pi_{1}|_{(\lambda , x)}$
 is tangent to $\frakh^* \times xH$,
on which  it is  isomorphic to the standard 
  Poisson (symplectic) structure on the cotangent bundle $T^* H$  (see, e.g.,
\cite{Gutt}).
Here $T^* H$ is identified with $\frakh^* \times H$  (hence 
with $\frakh^* \times xH$) via left translations.
It thus follows that $[\pi_{1} , \pi_{1}]=0$. Therefore
$$[\pi , \pi ]=2[\pi_{1} , \Vec{r(\lambda )}]+[\Vec{r(\lambda )}, \Vec{r(\lambda )}]. $$
Now
\be
&&[\pi_{1} , \  \Vec{r(\lambda )}]\\
&=&[\pi_{\frakh^* } , \  \Vec{r(\lambda )}]+\sum_{i} [\parr{}\wedge \vh{i}, \ 
\Vec{r(\lambda )}]\\
&=&[\pi_{\frakh^* } , \  \Vec{r(\lambda )}]+\sum_{i}
[\Vec{r(\lambda )},  \ \parr{}]\wedge \vh{i}- \sum_{i}
\parr{}\wedge [\Vec{r(\lambda )}, \  \vh{i}].
\ee
Hence $[\pi , \pi ]=I_{1}+I_{2}$, where
\be
I_{1}&=&2\sum_{i}
[\Vec{r(\lambda )}, \parr{}]\wedge \vh{i}+[\Vec{r(\lambda )}, \Vec{r(\lambda 
)}], \ \ \ \mbox{and}\\
I_{2}&=&2[\pi_{\frakh^* } , \Vec{r(\lambda )}]-2\sum_{i}
\parr{}\wedge [\Vec{r(\lambda )}, \vh{i}].
\ee
With respect
to the natural bigrading on $\wedge^{3}T(\frakh^* \times G)$, $I_{1}$ and $I_{2}$ correspond to the $(0, 3)$
and $(1, 2)$-terms   of $[\pi , \pi ]$, respectively.
It thus follows that $[\pi , \pi ]=0$ iff
$I_{1}=0$ and $I_{2}=0$.

It is simple to see that
$$I_{1}=-2\sum_{i}  \vh{i} \wedge \frac{\partial \Vec{r}}{ \partial\lambda^{i}}
+\Vec{[r(\lambda ), r(\lambda )]}.$$
Hence $I_{1}=0$ is equivalent to Equation (\ref{eq:triangular}).

To  find out the meaning of  $I_{2}=0$, let us write 
$\pi_{\frakh^* }=\half 
\sum_{ij}f_{ij}(\lambda )\frac{\partial }{\partial\lambda^{i}}
\wedge \frac{\partial }{\partial\lambda^{j}} \ (f_{ij}=-f_{ji})$. 
A simple computation yields that
\be
I_{2}&=&2 \sum_{i} \parr{}\wedge \sum_{j}f_{ij}(\lambda )\frac{\partial \Vec{r}}{ \partial\lambda^{j}} +2\sum_{i} \parr{}\wedge
\Vec{[h_{i}, r(\lambda )]}.
\ee
Thus $I_{2}=0$ is equivalent to
 $$[h_{i}, r(\lambda )]=-
\sum_{j}f_{ij}(\lambda )\frac{\partial  r(\lambda )
}{ \partial\lambda^{j}}= \left.\frac{d}{dt}\right|_{t=0} 
r(Ad^{*}_{\exp^{-1}{th_{i}}}\lambda ), \ \forall i,$$
which   exactly means that  $r$ is
 $H$-equivariant.  This concludes the proof. \qed
{\bf Remark.}\ Note that $M \ (=\frakh^* \times G)$
 admits a left $G$-action and a right $H$-action defined as follows:
$\forall (\lambda , x)\in \frakh^* \times G$,
\be
&&y \cdot (\lambda , x)=(\lambda , yx), \ \ \ \forall y\in G;\\
&&(\lambda , x)  \cdot y=(Ad_{y}^{*}\lambda , xy),  \ \ \ \forall y\in H.
\ee

It is  clear that the Poisson structure $\pi$
is invariant under both actions.
And, in short, we will say that $\pi$ is $G\times H$-invariant.

\begin{defi}
A classical triangular dynamical $r$-matrix $r: \frakh^* \lon
\wedge^2 \frakg$  is said to be
{\em non-degenerate} if the corresponding Poisson structure $\pi $  
on $M$ is non-degenerate, i.e., symplectic.
\end{defi}

In what follows, we will give  a geometric construction
of non-degenerate dynamical r-matrices. To this end, let us
first  recall a useful construction of a symplectic manifold
from a fat principal  bundle \cite{GLS, Weinstein80}. A principal
bundle $P(M, H)$ with  a connection  is called {\em fat} on
an open submanifold $U\subseteq \frakh^*$ if the
scalar-valued forms $<\lambda , \Omega >$ is 
non-degenerate  on each horizontal
space in $TP$ for $\lambda \in U$. Here $\Omega $ is
the curvature form, which is a tensorial 
form of type $Ad_H$ on $P$ (i.e., it is horizontal,
$\frakh$-valued, and $Ad_H$-equivariant).  

Given a  fat bundle $P(M, H)$ with  a connection, one
has a  decomposition of  the tangent bundle  $TP =Vert (P) \oplus
 Hor (P) $.
We may identify $Vert (P)$ with a trivial bundle with fiber
$\frakh$. Thus 
$$Vert^* P\cong  \frakh^* \times P.$$
On the other hand, $Vert^* P \cong Hor^{\perp} (P)\subset 
T^* P$. Thus, by pulling back the canonical symplectic
structure on $T^* P$, one can equip $Vert^* P$,
hence $ \frakh^* \times P$,  an $H$-invariant presymplectic structure,
where $H$ acts on $ \frakh^* \times P$ by $(\lambda , x)\cdot h
=(Ad_{h}^{*}\lambda , x\cdot h )$, $\forall h\in H$ and $(\lambda , x)
\in \frakh^* \times P$.
If $U\subseteq \frakh^*$ is an  open submanifold
on which $P(M, H)$ is fat,
 then we obtain an  $H$-invariant  symplectic manifold $U\times P$.
In fact, the presymplectic form $\omega $ can be described explicitly.
Note that $Vert^* P$ admits a natural fibration with 
$T^* H$ being the fibers,
 and the connection on  $P$ induces a connection on
this fiber bundle.
 In other words, $Vert^* P$ is a 
symplectic fibration in the sense  of  
Guillemin-Lerman-Sternberg \cite{GLS}.
At any point $(\lambda , x)\in \frakh^*\times P\cong 
Vert^* P$, the presymplectic form $\omega $ can be described
as  follows: it restricts to the canonical two-form
on the   fiber; the vertical subspace is $\omega$-orthogonal
to the horizontal subspace; and the horizontal subspace is
isomorphic to the horizontal subspace of $T_{x}P$ and the restriction
of $\omega$ to this subspace is the two form $-<\lambda , \Omega (x)>$
obtained by pairing the curvature form with $\lambda $
(see Examples 2.2-2.3 in \cite{GLS}). 
 
Now assume that
\begin{equation}
\label{eq:reductive}
\frakg =\frakh \oplus \frakm 
\end{equation}
is   a reductive decomposition
of a Lie algebra $\frakg$, i.e., $\frakh$ is a Lie subalgebra
and $\frakm$ is stable under the adjoint action of $\frakh$:
$[\frakh , \ \frakm ]\subset  \frakm$. By $G$, we denote a
Lie group with Lie algebra $\frakg$, and $H$  the  Lie subgroup
corresponding to $\frakh$. It is standard \cite{KN}
that the decomposition (\ref{eq:reductive})
 induces a  left $G$-invariant
connection on the principal bundle $G(G/H , H)$, where
the curvature is given by
\begin{equation}
\Omega (X, Y)=-[X, \ Y]_{\frakh}, \ \ \ \ \frakh-\mbox{component
of  } [X, \ Y]\in \frakg .
\end{equation}
Here $X$ and $Y$ are arbitrary left invariant vector fields 
on $G$ belonging to $\frakm$.

A reductive decomposition $\frakg =\frakh \oplus \frakm$
is said to be  {\em fat} if the corresponding principal bundle
$G(G/H , H)$ is fat on an open submanifold $U\subseteq \frakh^*$.
As a consequence, a fat decomposition $\frakg =\frakh \oplus \frakm$
gives rise to a $G\times H$-invariant symplectic structure
on $M= U\times G$, where the symplectic structure
is the restriction of the canonical symplectic
form on $T^* G$. In other words, $M$  is a symplectic submanifold of $T^*G$. 
Here the  embedding $U\times G \subseteq \frakh^* \times G
\lon \frakg^*  \times G \  (\cong T^*G )$ is 
 given by the natural inclusion $(\lambda , x )\lon (pr^* \lambda  , x)$, where
$pr: \frakg \lon \frakh$ is the projection along the decomposition
$\frakg =\frakh \oplus \frakm$.  Since the symplectic
structure $\omega$  on $U\times G$ is left invariant, in order to describe
$\omega$ explicitly,
it suffices to specify it  at a point $(\lambda , 1)$.
Now $T_{(\lambda , 1)}(U\times G )\cong \frakh^{*} \oplus \frakg
=\frakh^{*}\oplus \frakh \oplus \frakm$.
Under this identification, we have $\omega =\omega_{1}\oplus \omega_{2}$, where $\omega_{1}\in
\Omega^{2} (\frakh^{*}\oplus \frakh )$ is the canonical
symplectic two-form on $T^* H$ at the  point $(\lambda , 1)\in 
\frakh^{*}\times H \ (\cong T^* H)$,  and $\omega_{2} \in \Omega^{2} (\frakm )$
  is given by
$$\omega_{2}  (X, Y)=<\lambda ,  [X, \ Y]_{\frakh}>, \  \ \ \forall
X, Y\in \frakm. $$
Let $r(\lambda )\in \wedge^2 \frakm$ be the inverse of
$\omega_{2}$, which always exists for $\lambda
\in U$  since $\omega_{2}$ is
assumed to be non-degenerate on $U$. It thus follows that the Poisson structure
on $U\times G$ is  
$$\pi   = \pi_{ \frakh^* } + \sum_{i} \parr{}\wedge \vh{i} +\Vec{r(\lambda )}.$$
According to Proposition \ref{pro:dy},
 $r: U\lon \wedge^2 \frakm\subset  \wedge^2 \frakg$
is  a non-degenerate triangular dynamical r-matrix. Thus
we have proved

\begin{thm}
\label{thm:r}
Assume that $\frakg =\frakh \oplus \frakm$ is  a reductive decomposition
which is fat on an open submanifold  $U\subseteq \frakh^*$.
 Then the dual of the linear map
$\phi: \wedge^{2} \frakm \lon \frakh: \ (X, Y)\lon
[X, \ Y]_{\frakh}, \ \ \forall X, Y\in \frakm$
  defines  a non-degenerate triangular dynamical $r$-matrix
$r: U (\subseteq \frakh^* ) \lon \wedge^2  \frakm \subset \wedge^2 \frakg$,
$\forall \lambda \in U$.
Here $\frakm^*$ is identified  with $\frakm$ using the  non-degenerate bilinear
form $\phi^* (\lambda )\in \wedge^2  \frakm^*$.
\end{thm}

It is often more useful to express $r (\lambda )$
explicitly in terms of a basis.   To this end,
let us choose a  basis $\{e_{1}, \cdots , e_{m}\}$ of $\frakm$.
Let $a_{ij}(\lambda )=<\lambda , [e_{i}, e_{j}]_{\frakh}>,\ \ i, j=1,
\cdots , m$. By $(c_{ij}(\lambda ) )$ we denote the inverse of the matrix
 $(a_{ij}(\lambda )), \ \  \forall \lambda \in U$.
 Then one has 
\begin{equation}
\label{eq:ii}
r(\lambda )=\half \sum_{ij} c_{ij}(\lambda ) e_{i}\wedge e_{j}.
\end{equation}
{\bf Remark.}
\begin{enumerate}
\item  After the completion of the
first draft,  we learned  that a similar formula is also obtained
 independently by Etingof \cite{E}.
 Note that this dynamical r-matrix $r$ is   always   singular  at $0$.
To remove this singularity, one needs to make a shift of
the dynamical parameter $\lambda \to \lambda -\mu$.
\item It would be interesting to compare our formula with Theorem 3 in
\cite{ES}.
\end{enumerate}

We end this section with some examples.\\\\
{\bf Example 2.1}  Let $\frakg$ be a simple Lie algebra over $\complex$
and $\frakh$ a Cartan subalgebra. 
Let $$\frakg =\frakh \oplus \bigoplus_{\alpha\in
\Delta_{+}}(\frakg_{\alpha}\oplus \frakg_{-\alpha}) $$ be
the root  space decomposition, where $\Delta_{+}$ is the set of positive
roots with respect to $\frakh$.
Take $\frakm =\oplus_{\alpha\in
\Delta_{+}}(\frakg_{\alpha}\oplus \frakg_{-\alpha}) $.
Then $\frakg =\frakh \oplus \frakm$ is  clearly a
reductive decomposition. Let $e_{\alpha}\in
\frakg_{\alpha}$ and $e_{-\alpha}\in \frakg_{-\alpha}$ be 
dual  vectors with respect to the Killing form: $(e_{\alpha} , e_{-\alpha})
=1$.  For any $\lambda \in \frakh^*$,
set  $ a_{\alpha \beta}(\lambda )=< \lambda , [e_{\alpha} , 
e_{\beta}]_{\frakh}>,
\ \  \forall \alpha, \beta \in \Delta_{+}\cup (- \Delta_{+})$.
It is then clear that $a_{\alpha \beta}(\lambda )=0$, whenever
$\alpha +\beta \neq 0$; and 
\be
&&a_{\alpha,  -\alpha} (\lambda )\\
&=& <  \lambda , [e_{\alpha} ,e_{-\alpha}]_{\frakh}>\\
&=&(  \lambda , \alpha )( e_{\alpha} , e_{-\alpha})\\
&=&(  \lambda , \alpha ).
\ee 
Therefore, from Theorem \ref{thm:r} and Equation (\ref{eq:ii}), it follows that
$$r(\lambda)=-\sum_{\alpha \in \Delta_+} \frac{1}{( \lambda , \alpha )}
e_{\alpha} \wedge e_{-\alpha}$$
is a  non-degenerate triangular dynamical r-matrix, so we have
recovered this standard example in \cite{EV1}.\\\\\\\\
{\bf Example 2.2} 
As in the above example,
  let $\frakg$ be a simple Lie algebra over $\complex$ with a 
fixed Cartan subalgebra $\frakh$, and  $\frakl$  a  reductive
 Lie subalgebra  containing $\frakh$.
There is a subset $\Delta (\frakl )_{+}$
of $\Delta_{+}$ such that
$$\frakl =\frakh \oplus \bigoplus_{\alpha\in
\Delta (\frakl )_{+}}(\frakg_{\alpha}\oplus \frakg_{-\alpha}). $$ 
Let $\overline{\Delta}_{+}={\Delta}_{+}- \Delta (\frakl )_{+}, \ \
\Delta (\frakl )=\Delta (\frakl )_{+}\cup (-\Delta (\frakl )_{+})$,
and $\overline{\Delta}= \overline{\Delta}_{+}\cup (-\overline{\Delta}_{+})$, and
denote by $\frakm$ the subspace of $\frakg$:
$$\frakm =\sum_{\alpha\in
\overline{\Delta}_{+} } (\frakg_{\alpha}\oplus \frakg_{-\alpha}).$$
It is simple to see that  $\frakg=\frakl \oplus \frakm$ is 
indeed a fat reductive decomposition, and therefore 
induces a non-degenerate triangular dynamical r-matrix $r: \frakl^*
\lon \wedge^2 \frakg$. To  describe  $r$ explicitly, we note that
 the dual space $\frakl^*$ admits
a natural decomposition
$$\frakl^* =\frakh^* \oplus \bigoplus_{\alpha\in
\Delta (\frakl )_{+}}(\frakg_{\alpha}^{*}\oplus \frakg_{-\alpha}^{*}). $$
Hence  any element $\mu \in \frakl^* $ can be written 
as $\mu =\lambda \oplus \oplus_{\alpha \in \Delta (\frakl ) }\xi_{\alpha }$, where
$\lambda \in \frakh^*$ and $\xi_{\alpha }\in \frakg_{\alpha}^{*}$.
Let $a_{\alpha \beta }(\mu )=<\mu ,  [e_{\alpha} , e_{\beta}]_{\frakl}>,
\ \ \forall  \alpha, \beta \in \overline{\Delta}$.
It  is easy to see that
\begin{equation}
a_{\alpha \beta }(\mu )= \left\{
\begin{array}{ll}
(\lambda , \alpha ), &\mbox{if } \alpha +\beta =0;\\
<\xi_{\gamma}, [e_{\alpha} , e_{\beta}]>, & \mbox{if } \alpha +\beta  =
\gamma \in \Delta (\frakl );\\
0, &\mbox{otherwise.}  
\end{array}
\right.
\end{equation}
By $(c_{\alpha \beta} (\mu ))$, we denote the inverse matrix of
$(a_{\alpha \beta }(\mu ))$. According to Equation (\ref{eq:ii}),
$$r (\mu )=\half \sum_{\alpha , \beta \in \overline{\Delta}}
c_{\alpha \beta} (\mu ) e_{\alpha }\wedge e_{\beta}$$
is a   non-degenerate triangular dynamical r-matrix
over $\frakl^*$. In particular, if $\mu =\lambda \in \frakh^*$,
it follows immediately that 
\begin{equation}
\label{eq:rr}
r(\lambda)=-\sum_{\alpha \in \overline{\Delta}_+} \frac{1}{( \lambda , \alpha )}
e_{\alpha} \wedge e_{-\alpha} .
\end{equation}
Equation (\ref{eq:rr}) was first obtained by 
Etingof-Varchenko in \cite{EV1}. \\\\\\

The following example was pointed out to us by  D. Vogan.\\\\
{\bf Example 2.3}
Let $\frakg=\reals^{m+n}\oplus \reals^{m+n} \oplus \reals$ be
a $2(m +n)+1$ dimensional Heisenberg  Lie algebra and 
$\frakh =\reals^{n}\oplus \reals^{n} \oplus \reals$ 
its standard Heisenberg  Lie subalgebra. 
By $\{p_{i}, q_{i}, c \}, \ i=1, \cdots , n+m$, we denote
the standard generators of $\frakg$ and
 $\{p_{m+i}, q_{m+i}, c \}, \ i=1, \cdots , n$,
the generators of $\frakh$.
Let $\frakm$ be the subspace of $\frakg$ generated by
$\{p_{i}, q_{i}\},  \ i=1, \cdots , m$. It is 
then clear that $\frakg =\frakh \oplus \frakm$ 
is a reductive decomposition.
Let $\{p_{i}^* , q_{i}^* , c^* \}, \ i=1, \cdots , n+m$, 
be the dual basis corresponding to the standard generators of $\frakg$.
For any $\lambda \in \frakh^*$, write
$\lambda =\sum_{i=1}^{n} (a_{i} p_{m+i}^* +b_{i}q_{m+i}^{*})
+x c^*$. This induces    a coordinate system on 
$\frakh^*$, and therefore a  function on $\frakh^*$ can be identified
with a function with variables  $(a_{i}, b_{i}, x)$.
It is  clear that
\be
&&\omega  (p_{i} , q_{j} )(\lambda )=<\lambda , [p_{i} , q_{j}]_{\frakh}>=
x \delta_{ij};\\
&&\omega (p_{i}, p_{j})=\omega (q_{i}, q_{j})=0, \ \ \ \forall i,\ j=1, \cdots , m. 
\ee
It thus follows that
$$r(a_{i}, b_{i}, x)=-\frac{1}{x}\sum_{i=1}^{m} p_{i}\wedge q_{i}:
\ \ \frakh^* \lon \wedge^2 \frakg $$
is a  non-degenerate triangular dynamical r-matrix.
 
\section{Compatible star products}
From Proposition \ref{pro:dy},  we know that a triangular dynamical r-matrix
$r: \frakh^{*} \lon \wedge^{2} \frakg$ is equivalent
to a special type of  Poisson structures on $\frakh^* \times G$.
It is thus very natural to expect that  quantization
of $r$ can be derived from  a certain special type of
star-products on $\frakh^* \times G$.
It is simple to see that the Poisson brackets on $C^{\infty}(\frakh^* \times G)$
can be described as follows:
\begin{enumerate}
\item for any $f, g\in
C^{\infty}(\frakh^* )$,
$\{f, g \}=\{f, g \}_{\pi_{\frakh^* }}$;
\item for any $f \in C^{\infty}(\frakh^* )$ and $g\in C^{\infty}(G)$,
$\{f, g \}=\sum_{i} (\frac{\partial f}{\partial \lambda^{i}})(\vh{i}g )$;
\item for any $f, g \in C^{\infty}(G)$, $\{f, g \}=
\Vec{r(\lambda )}(f, g)$.
\end{enumerate}
These  Poisson bracket relations naturally  motivate  the
following:

\begin{defi}
\label{def:prefer}
A star product $\hstar$
 on $M=\frakh^* \times G$ is called a \prefer star product  if
\begin{enumerate}
\item for any $f, g\in C^{\infty}(\frakh^* )$,
\begin{equation}
\label{eq:1}
f(\lambda )*_{\hbar}g (\lambda )=f(\lambda ) \sstar g(\lambda );
\end{equation}
\item 
for any $f(x) \in  C^{\infty}(G)$ and $g(\lambda )\in C^{\infty}(\frakh^* )$,
\begin{equation}
\label{eq:2}
f(x)*_{\hbar}g (\lambda )=f(x)g (\lambda );
\end{equation}
\item
  for any $f (\lambda )\in C^{\infty}(\frakh^* )$
and $  g  (x)\in  C^{\infty}(G)$,
\begin{equation}
\label{eq:3}
f(\lambda )*_{\hbar}g (x)
=\sum_{k=0}^{\infty} \frac{ \hbar^k }{k!}
\frac{\partial^{k} f}{\partial \lambda^{i_{1}}\cdots \partial \lambda^{i_{k}}}
\vh{i_{1}}\cdots \vh{i_{k}}g;
\end{equation}
\item 
for any $f(x), \ g(x)\in  C^{\infty}(G)$,
\begin{equation}
\label{eq:4}
f(x)*_{\hbar}g(x)= \Vec{F(\lambda )}(f, g),
\end{equation}
where $F(\lambda )$ is a smooth function
 $F: \frakh^* \lon U\frakg \ot U\frakg \flb\hbar \frb$
such that $F= 1 +\hbar F_{1} + O(\hbar^2 )$.
\end{enumerate}
\end{defi}
Here $\sstar $ denotes the standard PBW-star product on $\frakh^*$
quantizing the canonical   Lie-Poisson
 structure   (see \cite{AW}), whose definition is recalled below. 
Let $\frakh_{\hbar}=\frakh \flb \hbar \frb$ be   a  Lie
algebra with the Lie bracket
$[X, Y]_\hbar =\hbar [X, Y]$, $\forall X, Y\in \frakh  \flb \hbar \frb$,
 and 
$$\sigma :  S(\frakh) \flb \hbar  \frb\cong U\frakh_{\hbar}$$
be the Poincar\'e-Birkhoff-Witt map, which is a vector space isomorphism.
Thus the  multiplication on $U\frakh_{\hbar}$ induces 
 a  multiplication
on $S(\frakh) \flb \hbar \frb\  (\cong \mbox{Pol}(\frakh^* )\flb \hbar \frb )$,
 hence on $C^{\infty} (\frakh^* )\flb \hbar \frb$, which
is  denoted by $\sstar $. It is easy to check that $\sstar $  satisfies
$$f\sstar  g= fg+\half \hbar \{f, g\}_{\pi_\frakh^* }+\sum_{k\geq 0}\hbar^{k}
B_{k}(f, g)+\cdots , \ \ \forall f, g\in C^{\infty} (\frakh^* ),$$
where $B_{k}$'s are bidifferential operators. In other words,
$\sstar $ is indeed a star product on $\frakh^*$, which is called the
 PBW-star product.

The following proposition is quite obvious.

\begin{pro}
\label{pro:prefer}
The classical limit of a \prefer star product is
the Poisson 
structure  $\pi = \pi_{\frakh^* } + \sum_{i}  \parr{}\wedge \vh{i}
 +\Vec{r(\lambda )}$,
where $r(\lambda )=F_{12}(\lambda )-F_{21}(\lambda )$.
\end{pro}

Below we will study some important properties of \prefer star products.

\begin{pro}
A \prefer star product is always  invariant under the left $G$-action.
It is right $H$-invariant iff $F: \frakh^* \lon U\frakg \ot U\frakg \flb\hbar \frb$
is $H$-equivariant, where $H$ acts on $\frakh^*$ by the coadjoint action
and on $U\frakg \ot U\frakg$ by the adjoint action.
\end{pro}
\pf  First of all, note that Equations (\ref{eq:1}-\ref{eq:4})
completely determine a star product.
It is clear, from these equations, that $\hstar$ is left $G$-invariant.

As for the right $H$-action, it is obvious from 
 Equation (\ref{eq:2}) that $*_\hbar  $  is invariant for 
$f(x)*_{\hbar}g (\lambda )$.
It is standard that  $\sstar $ is  invariant under the coadjoint action,
so it follows from  Equation (\ref{eq:1}) that
 $f(\lambda )*_{\hbar}g(\lambda )$
  is also $H$-invariant.

For any $h\in \frakh$, $g(x)\in C^{\infty}(G)$ and any fixed 
$y\in H$,

\be
\vh{} (R_{y}^{*}g)(x)&=&(L_{x}h)(R_{y}^{*}g)\\
&=&(R_{y}L_{x}h)(g )\\
&=&(L_{xy} Ad_{y^{-1}}h)(g )\\
&=&(\Vec{Ad_{y^{-1}}h}g)(xy)\\
&=&[R_{y}^{*}(\Vec{Ad_{y^{-1}}h}g)](x).
\ee
Thus it follows that
\begin{equation}
\label{eq:5}
\vh{i_{1}}\cdots \vh{i_{k}} (R_{y}^{*}  g)
=R_{y}^{*} (\vhh{i_{1}}\cdots \vhh{i_{k}} g), 
\end{equation}
where $h_{i}'=Ad_{y^{-1}}h_{i}$, $i=1,\ \cdots , n$.
Let $\xi_{i}'=Ad_{y}^{*}\xi_{i}$, $i=1,\ \cdots , n$.
Then $\{\xi_{1}', \cdots , \xi_{l}'\}$ is a  dual basis
for $\{h_{1}', \cdots , h_{l}'\}$.
Let $(\lambda^{'1}, \cdots , \lambda^{'l})$ be  its corresponding  
 induced coordinates  on $\frakh^*$. 
Then
\be
\frac{\partial}{\partial \lambda^{i}}((Ad_{y}^{*})^{*}f)&=&
\left.\frac{d}{dt}\right|_{t=0}  ((Ad_{y}^{*})^{*}f)(\lambda +t\xi_{i} )\\
&=&\left.\frac{d}{dt}\right|_{t=0}  f(Ad_{y}^{*}\lambda +tAd_{y}^{*}\xi_{i} )\\
&=&\left.\frac{d}{dt}\right|_{t=0}  f(Ad_{y}^{*}\lambda +t\xi_{i}' )\\
&=&\frac{\partial f}{\partial \lambda^{'i}} (Ad_{y}^{*}\lambda )\\
&=&(Ad_{y}^{*})^{*}\frac{\partial f}{\partial \lambda^{'i}} .
\ee
Hence 
\begin{equation}
\label{eq:6}
\frac{\partial^{k} [(Ad_{y}^{*})^{*}f]}{\partial \lambda^{i_{1}}\cdots \partial \lambda^{i_{k}}}
=(Ad_{y}^{*})^{*}[\frac{\partial^{k} f}{\partial \lambda^{'i_{1}}\cdots \partial
\lambda^{'i_{k}}}].
\end{equation}
 From Equation (\ref{eq:3}), it  follows that
for any $f (\lambda )\in C^{\infty}(\frakh^* )$
and $  g  (x)\in  C^{\infty}(G)$,
\be
(R_{y}^* f)(\lambda )*_{\hbar} (R_{y}^* g) (x)
&=&\sum_{k=0}^{\infty} \frac{ \hbar^k }{k!}
\frac{\partial^{k} 
[(Ad_{y}^{*})^{*}f] }{\partial \lambda^{i_{1}}\cdots \partial \lambda^{i_{k}}}
\vh{i_{1}}\cdots \vh{i_{k}} (R_{y}^* g) \ \ \mbox{(by 
Equations (\ref{eq:5}-\ref{eq:6}))}\\
&=&\sum_{k=0}^{\infty} \frac{ \hbar^k }{k!}
  (Ad_{y}^{*})^{*}[\frac{\partial^{k} f }{\partial \lambda^{'i_{1}}
\cdots \partial \lambda^{'i_{k}}} ]
R_{y}^* [\Vec{h'_{i_{1}}}\cdots \Vec{h'_{i_{k}}}  g]\\
&=&R_{y}^* (f(\lambda )*_{\hbar}  g  (x) ).
\ee
I.e., $f(\lambda )*_{\hbar}  g  (x)$ is also right $H$-invariant.

Finally,  $\forall f(x), \ g(x)\in  C^{\infty}(G)$,
\be
&&(R_{y}^{*}(f*_{\hbar}g))(\lambda , x)\\
&=&(f*_{\hbar}g )(Ad_{y}^{*}\lambda , xy)\\
&=&\Vec{F(Ad_{y}^{*}\lambda )}(f, g)(xy)\\
&=&[L_{xy}(F(Ad_{y}^{*}\lambda ))](f, g ).
\ee

On the other hand,
\be
&&(R_{y}^{*}f*_{\hbar}R_{y}^{*}g )(\lambda , x)\\
&=&\Vec{F(\lambda )}(R_{y}^{*}f , R_{y}^{*}g )(x)\\
&=&(L_{x}F(\lambda )) (R_{y}^{*}f , R_{y}^{*}g )\\
&=&(R_{y}L_{x}F(\lambda ) ) (f, g).
\ee
Therefore $R_{y}^{*}(f*_{\hbar}g )=R_{y}^{*}f*_{\hbar}R_{y}^{*}g$
iff $L_{xy}(F(Ad_{y}^{*}\lambda ))=R_{y}L_{x}F(\lambda )$.  The latter
is equivalent to that $F(Ad_{y}^{*}\lambda )=Ad_{y^{-1}}F(\lambda )$,
or $F$ is $H$-equivariant. This concludes the proof. \qed

In order to give an explicit formula for $*_{\hbar}$, let us
write 
\begin{equation}
F(\lambda )=\sum a_{\alpha \beta} (\lambda ) U_{ \alpha }\ot U_{\beta},
\end{equation}
where $a_{\alpha \beta} (\lambda )\in C^{\infty}(\frakh^* )\flb\hbar \frb$
and $U_{ \alpha }\ot U_{\beta}\in U\frakg \ot U\frakg $.
Using this expression, indeed one  can  describe $*_{\hbar}$
explicitly.


\begin{thm}
\label{thm:formula}
Given a \prefer star product $*_{\hbar} $ as in Definition \ref{def:prefer},
  for any $f(\lambda , x), \ 
g(\lambda , x)\in C^{\infty}(\frakh^* \times G)\flb\hbar \frb$,
\begin{equation}
\label{eq:explicit}
f(\lambda , x)*_{\hbar} g(\lambda , x)=
\sum_{\alpha \beta} \sum_{k=0}^{\infty} \frac{ \hbar^k }{k!}
a_{\alpha \beta} (\lambda )\sstar \Vec{U_{ \alpha }}
 \frac{\partial^{k} f}{\partial \lambda^{i_{1}}\cdots \partial \lambda^{i_{k}}}
\sstar \Vec{U_{\beta}} \vh{i_{1}}\cdots \vh{i_{k}}g .
\end{equation}
\end{thm}

We need a couple of lemmas first.

\begin{lem}
Under the same hypothesis as in Theorem \ref{thm:formula},
\begin{enumerate}
\item for any $f(\lambda , x) \in  C^{\infty}(\frakh^* \times G )$
 and $g(\lambda )\in   C^{\infty}(\frakh^* )$,
\begin{equation}
\label{eq:7}
f(\lambda , x )*_{\hbar}g (\lambda )=f(\lambda , x ) \sstar g(\lambda );
\end{equation}
\item for any $f(x) \in  C^{\infty}(G)$ and $g(\lambda , x) \in  
C^{\infty} (\frakh^* \times G)$,
\begin{equation}
\label{eq:8}
f(x)*_{\hbar} g(\lambda , x)=\sum_{\alpha \beta}a_{\alpha \beta} (\lambda )
\sstar \Vec{U_{ \alpha }}f (x) \Vec{U_{\beta}}g (\lambda , x);
\end{equation}
\item for any $f(\lambda , x) \in C^{\infty} (\frakh^* \times G)$
and $g(x)\in C^{\infty}(G)$,
\begin{equation}
\label{eq:9}
f(\lambda , x )*_{\hbar}g (x)=
\sum_{\alpha \beta} \sum_{k=0}^{\infty} \frac{ \hbar^k }{k!}
a_{\alpha \beta} (\lambda )\sstar \Vec{U_{ \alpha }}
 \frac{\partial^{k} f (\lambda , x )}{\partial
 \lambda^{i_{1}}\cdots \partial \lambda^{i_{k}}} 
\Vec{U_{\beta}} \vh{i_{1}}\cdots \vh{i_{k}}g (x).
\end{equation}
\end{enumerate}
\end{lem}
\pf  (i). It suffices to show this identity for 
$f(\lambda , x)= f_{1}(x )f_{2}( \lambda)$,
$\forall f_{1}(x )\in
C^{\infty}(G )$ and $f_{2}(\lambda )\in C^{\infty}(\frakh^* )$.
Now
\be
&&f(\lambda , x)*_{\hbar}g (\lambda )\\
&=&(f_{1}(x )f_{2}( \lambda))*_{\hbar}g (\lambda )\ \ \ 
 \mbox{(by Equation (\ref{eq:2}))}\\
&=&(f_{1}(x )*_{\hbar}f_{2}( \lambda))*_{\hbar}g (\lambda )\\
&=&f_{1}(x )*_{\hbar}(f_{2}( \lambda)*_{\hbar}g (\lambda )) \ \ \mbox{(by
Equations (\ref{eq:1}-\ref{eq:2}))}\\
&=&f_{1}(x )(f_{2}( \lambda)\sstar g (\lambda ))\\
&=&(f_{1}(x )f_{2}( \lambda))\sstar g (\lambda )\\
&=&f(\lambda , x)\sstar g (\lambda ).
\ee

(ii). Similarly, we may assume that $g(\lambda , x)
=g_{1}(x )g_{2}( \lambda)$, $\forall g_{1}(x )\in
C^{\infty}(G )$ and $g_{2}(\lambda )\in C^{\infty}(\frakh^* )$. 
Then,
\be
&&f(x)*_{\hbar}g(\lambda , x)\\
&=&f(x)*_{\hbar} (g_{1}(x )g_{2}( \lambda))\\
&=&f(x)*_{\hbar}(g_{1}(x )*_{\hbar}g_{2}( \lambda))\\
&=&(f(x)*_{\hbar}g_{1}(x ))*_{\hbar}g_{2}( \lambda ) \ \ \ \mbox{(by Equation
(\ref{eq:4}))}\\
&=&\sum_{\alpha \beta} [a_{\alpha \beta} (\lambda ) (\Vec{U_{ \alpha }}f(x))
 (\Vec{U_{\beta}} g_{1}(x ))]\sstar g_{2}( \lambda )\\
&=&\sum_{\alpha \beta}a_{\alpha \beta} (\lambda ) \sstar 
\Vec{U_{ \alpha }}f (x)  \Vec{U_{\beta}}g (\lambda , x).
\ee

(iii). Assume that $f(\lambda , x)= f_{1}(x )f_{2}( \lambda), \ \forall f_{1}(x )\in C^{\infty}(G )$
 and $f_{2}(\lambda )\in C^{\infty}(\frakh^* )$.
Then
\be
&&f(\lambda , x )*_{\hbar}g (x)\\
&=&(f_{1}(x )f_{2}( \lambda))*_{\hbar}g (x)\\
&=&(f_{1}(x )*_{\hbar}f_{2}( \lambda))*_{\hbar}g (x)\\
&=&f_{1}(x )*_{\hbar}(f_{2}( \lambda)*_{\hbar}g (x)) \ \ \ \ \mbox{(using 
Equation (\ref{eq:8}))}\\
&=&\sum_{\alpha \beta}a_{\alpha \beta} (\lambda ) \sstar
\Vec{U_{ \alpha }}f_{1}(x ) \Vec{U_{\beta}}(f_{2}( \lambda)*_{\hbar}g (x))
\\
&=&\sum_{\alpha \beta} \sum_{k=0}^{\infty} \frac{ \hbar^k }{k!}
a_{\alpha \beta} (\lambda )\sstar [\Vec{U_{ \alpha }} f_{1}(x )
\Vec{U_{\beta}}(
 \frac{\partial^{k} f_{2}( \lambda)}{\partial \lambda^{i_{1}}\cdots \partial
 \lambda^{i_{k}}} \vh{i_{1}}\cdots \vh{i_{k}}g(x)  )]\\
&=&\sum_{\alpha \beta} \sum_{k=0}^{\infty} \frac{ \hbar^k }{k!}
a_{\alpha \beta} (\lambda )\sstar [\Vec{U_{ \alpha }} f_{1}(x )
\frac{\partial^{k} f_{2}( \lambda)}{\partial \lambda^{i_{1}}\cdots \partial 
\lambda^{i_{k}}}
\Vec{U_{\beta}} \vh{i_{1}}\cdots \vh{i_{k}}g(x)]\\
&=&\sum_{\alpha \beta} \sum_{k=0}^{\infty} \frac{ \hbar^k }{k!}
a_{\alpha \beta} (\lambda )\sstar [\Vec{U_{ \alpha }}
\frac{\partial^{k} (f_{1}(x )f_{2}( \lambda ))}{\partial
 \lambda^{i_{1}}\cdots \partial \lambda^{i_{k}}}
\Vec{U_{\beta}} \vh{i_{1}}\cdots \vh{i_{k}}g(x)]\\
&=&\sum_{\alpha \beta} \sum_{k=0}^{\infty} \frac{ \hbar^k }{k!}
a_{\alpha \beta} (\lambda )\sstar \Vec{U_{ \alpha }}
 \frac{\partial^{k} f (\lambda , x)}{\partial \lambda^{i_{1}}\cdots \partial \lambda^{i_{k}}}
 \Vec{U_{\beta}} \vh{i_{1}}\cdots \vh{i_{k}}g(x).
\ee
This concludes the proof of the lemma. \qed

Now we are ready to prove the main result of this section.\\\\\\
{\bf Proof of Theorem \ref{thm:formula}} Again,  we 
may assume that $g(\lambda , x)
=g_{1}(x )g_{2}( \lambda)$, $\forall g_{1}(x )\in
C^{\infty}(G )$ and $g_{2}(\lambda )\in C^{\infty}(\frakh^* )$.
Then
\be
&&f(\lambda , x)*_{\hbar} g(\lambda , x)\\
&=&f(\lambda , x)*_{\hbar} (g_{1}(x )g_{2}( \lambda))\\
&=&f(\lambda , x)*_{\hbar} (g_{1}(x )*_{\hbar}g_{2}( \lambda))\\
&=&(f(\lambda , x)*_{\hbar}g_{1}(x ))*_{\hbar}g_{2}( \lambda)\ \ \ \ 
\mbox{(by Equation (\ref{eq:7}))}\\
&=&(f(\lambda , x)*_{\hbar}g_{1}(x ))\sstar g_{2}( \lambda )
\ \ \ \mbox{(by Equation (\ref{eq:9}))}\\
&=& \sum_{\alpha \beta} \sum_{k=0}^{\infty} \frac{ \hbar^k }{k!}
[a_{\alpha \beta} (\lambda )\sstar \Vec{U_{ \alpha }}
 \frac{\partial^{k} f (\lambda , x)}{\partial \lambda^{i_{1}}\cdots \partial \lambda^{i_{k}}}
 \Vec{U_{\beta}} \vh{i_{1}}\cdots \vh{i_{k}}g_{1}(x)]\sstar g_{2}( \lambda)\\
&=& \sum_{\alpha \beta} \sum_{k=0}^{\infty} \frac{ \hbar^k }{k!}
a_{\alpha \beta} (\lambda )\sstar \Vec{U_{ \alpha }}
 \frac{\partial^{k} f (\lambda , x)}
{\partial \lambda^{i_{1}}\cdots \partial \lambda^{i_{k}}}
\sstar \Vec{U_{\beta}} \vh{i_{1}}\cdots \vh{i_{k}}(g_{1}(x) g_{2}( \lambda))\\
&=&\sum_{\alpha \beta} \sum_{k=0}^{\infty} \frac{ \hbar^k }{k!}
a_{\alpha \beta} (\lambda )\sstar \Vec{U_{ \alpha }}
 \frac{\partial^{k} f ( \lambda , x)}{\partial \lambda^{i_{1}}\cdots \partial \lambda^{i_{k}}}
\sstar \Vec{U_{\beta}} \vh{i_{1}}\cdots \vh{i_{k}}g ( \lambda , x).
\ee
\qed

As a  consequence of Theorem \ref{thm:formula}, we will  see that 
 if a  function $F(\lambda ):
 \frakh^* \lon U\frakg \otimes U\frakg \flb \hbar \frb$
defines  a  \prefer star product, it must  satisfy
a ``twisted-cocycle" type condition.
To describe this condition explicitly, we need to  introduce some
notations.

For any    $f(\lambda )\in C^{\infty}(\frakh^* )$, define
$f(\lambda +\hbar h)\in C^{\infty}(\frakh^* )\ot U\frakh \flb\hbar \frb$ by
\begin{eqnarray}
f (\lambda +\hbar h ) &=&f(\lambda  )\ot 1+\hbar \sum_{i}
\parr{f}\ot h_{i}       +\frac{1}{2!}\hbar^{2}\sum_{i_{1}i_{2}}
\frac{\partial^{2}f}{\partial \lambda_{i_{1}}\partial \lambda_{i_{2}}}\ot
h_{i_{1}}h_{i_{2}} \nonumber\\
&& \ \ \ \ +\cdots
+\frac{\hbar^k}{k!}\sum \frac{\partial^{k} f}{\partial  \lambda_{i_{1}}
\cdots \partial \lambda_{i_{k}}} \ot h_{i_{1}}\cdots h_{i_{k}} +\cdots.
\end{eqnarray}

The correspondence  $C^{\infty}(\frakh^* )
\lon C^{\infty}(\frakh^* )\ot U\frakh \flb\hbar \frb:\  f(\lambda )\lon
 f (\lambda +\hbar h )$
extends naturally to a  linear map
 from $C^{\infty}(\frakh^* ) \ot U\frakg \ot U\frakg \flb\hbar \frb$
to $C^{\infty}(\frakh^* )\ot U\frakh \ot U\frakg \ot U\frakg \flb\hbar \frb
\subseteq C^{\infty}(\frakh^* )\ot U\frakg \ot U\frakg \ot U\frakg \flb\hbar \frb$,
which is denoted by $F(\lambda )\lon F_{23}( \lambda +\hbar h^{(1)})$.
More explicitly, assume that $F(\lambda  )=\sum_{\alpha \beta}
f_{\alpha \beta}( \lambda )U_{\alpha} \ot U_{\beta}$, where
$f_{\alpha \beta}( \lambda )\in C^{\infty}(\frakh^* )\flb\hbar \frb$
and $U_{\alpha} \ot U_{\beta}\in U\frakg \ot U\frakg$.
Then
\begin{equation}
F_{23}( \lambda +\hbar h^{(1)})=\sum_{\alpha \beta}f_{\alpha \beta}
(\lambda +\hbar h )\ot U_{\alpha} \ot U_{\beta}.
\end{equation}

By a suitable permutation, one  may define $F_{12}( \lambda +\hbar h^{(3)})$ 
and $F_{13}( \lambda +\hbar h^{(2)})$ similarly. Note  that
 $U\frakg$ is a   Hopf algebra. By  
$\Delta : U\frakg \lon U\frakg \ot U\frakg$ and
 $\epsilon : U\frakg \lon \reals$, we denote its
co-multiplication and  co-unit, respectively.
Then $\Delta $  naturally extends
to a map $C^{\infty}(\frakh^* )\ot U\frakg \flb\hbar \frb
 \lon C^{\infty}(\frakh^* ) \ot U\frakg \ot U\frakg \flb\hbar \frb$, which will
be denoted by the same symbol.
 
\begin{cor}
\label{cor:shift}
Assume that  $F: \frakh^* \lon U\frakg \ot U\frakg \flb\hbar \frb$
defines  a \prefer star product $*_{\hbar}$ as in Definition
\ref{def:prefer}. Then
\begin{eqnarray}
&&(\Delta  \ot  id )F(\lambda ) \sstar  F_{12} (\lambda +\hbar h^{(3)})
   =  (id \ot  \Delta )  F (\lambda ) \sstar F_{23}(\lambda  );
\label{eq:shifted0}\\
&&\label{eq:co0}
   (\epsilon \ot id) F(\lambda )  =  1; \ \
(id \ot \epsilon ) F (\lambda ) =  1.
\end{eqnarray}
\end{cor}
\pf Equation (\ref{eq:co0}) follows from the fact that
$1*_{\hbar}f(x)=f(x)*_{\hbar}1=f(x), \ \forall f(x)\in C^{\infty}(G)$.

As for Equation (\ref{eq:shifted0}),
 note that for any $f_{1}(x), f_{2}(x)$ and $ f_{3}(x)\in C^{\infty}(G)$,
according to Equation (\ref{eq:9}), we have
\be
&&(f_{1}(x)*_{\hbar}f_{2}(x)) *_{\hbar}f_{3}(x)\\
&=&\sum_{\alpha \beta} \sum_{k=0}^{\infty} \frac{ \hbar^k }{k!}
a_{\alpha \beta} (\lambda )\sstar \Vec{U_{ \alpha }}
 \frac{\partial^{k} (f_{1}(x)*_{\hbar}f_{2}(x))}{\partial \lambda^{i_{1}}\cdots \partial \lambda^{i_{k}}}
 \Vec{U_{\beta}} \vh{i_{1}}\cdots \vh{i_{k}}f_{3}(x).
\ee

Now
\be
&&(\Delta  \ot  id )F(\lambda ) \sstar  F_{12} (\lambda +\hbar h^{(3)})\\
&=& \sum_{k=0}^{\infty} \frac{ \hbar^k }{k!}
 (\Delta  \ot  id )F(\lambda ) \sstar 
 (\frac{\partial^{k} F}{\partial \lambda^{i_{1}}\cdots \partial \lambda^{i_{k}}}
 \ot  h_{i_{1}}\cdots h_{i_{k}})\\
&=&\sum_{\alpha \beta} \sum_{k=0}^{\infty} \frac{ \hbar^k }{k!}
a_{\alpha \beta} (\lambda )\sstar \Delta U_{ \alpha } 
\frac{\partial^{k} F}{\partial \lambda^{i_{1}}\cdots \partial \lambda^{i_{k}}}
\ot U_{\beta}h_{i_{1}}\cdots h_{i_{k}} .
\ee
It thus follows that
\be
&&\Vec{(\Delta  \ot  id )F(\lambda ) \sstar  F_{12} (\lambda +\hbar h^{(3)})}
(f_{1}(x),\ f_{2}(x) ,\ f_{3}(x))\\
&=&\sum_{\alpha \beta} \sum_{k=0}^{\infty} \frac{ \hbar^k }{k!}
a_{\alpha \beta} (\lambda )\sstar
\Vec{U_{ \alpha } }
(\frac{\partial^{k} \Vec{F(\lambda )} (f_{1}(x), f_{2}(x))}{\partial \lambda^{i_{1}}\cdots \partial \lambda^{i_{k}}} )
\Vec{U_{\beta}} \vh{i_{1}}\cdots \vh{i_{k}}f_{3}(x)\\
&=&\sum_{\alpha \beta} \sum_{k=0}^{\infty} \frac{ \hbar^k }{k!}
a_{\alpha \beta} (\lambda )\sstar \Vec{U_{ \alpha }}
 \frac{\partial^{k} (f_{1}(x)*_{\hbar}f_{2}(x))}{\partial \lambda^{i_{1}}\cdots
\partial \lambda^{i_{k}}}
 \Vec{U_{\beta}} \vh{i_{1}}\cdots \vh{i_{k}}f_{3}(x)\\
&=&(f_{1}(x)*_{\hbar}f_{2}(x)) *_{\hbar}f_{3}(x).
\ee

On the other hand,
\be
&&f_{1}(x)*_{\hbar}(f_{2}(x) *_{\hbar} f_{3}(x))\\
&=&f_{1}(x)*_{\hbar} \Vec{F(\lambda )} (f_{2}(x), f_{3}(x)) \ \ \mbox{(by Equation (\ref{eq:8}))} \\
&=&\sum_{\alpha \beta}a_{\alpha \beta} (\lambda )
\sstar \Vec{U_{ \alpha }}f_{1}(x)  \Vec{U_{\beta}} (\Vec{F(\lambda )} (f_{2}(x), f_{3}(x)) )\\
&=&\Vec{(id \ot  \Delta )  F (\lambda ) \sstar F_{23}(\lambda  )}
(f_{1}(x),\ f_{2}(x) ,\ f_{3}(x)).
\ee
Now Equation (\ref{eq:shifted0})  follows from
 the associativity of $*_{\hbar}$.\qed

To end this section,
as a special case, let us consider $M=\frakh^* \times H\cong T^* H$, which is
equipped with the canonical cotangent   symplectic structure.
The following proposition  describes an explicit formula for
a \prefer star-product on it.


\begin{pro}
\label{pro:H}
 For any $f(\lambda , x), \
g(\lambda , x)\in C^{\infty}( \frakh^* \times H )\flb\hbar \frb$,
the following equation
\begin{equation}
\label{eq:19}
f(\lambda , x)*_{\hbar} g(\lambda , x)=
 \sum_{k=0}^{\infty} \frac{ \hbar^k }{k!}
 \frac{\partial^{k} f}{\partial \lambda^{i_{1}}\cdots \partial \lambda^{i_{k}}}
\sstar \vh{i_{1}}\cdots \vh{i_{k}}g
\end{equation}
defines a \prefer star product on $M=\frakh^* \times H\cong T^* H$,
which is in fact  a deformation quantization of  its
canonical cotangent   symplectic structure.
\end{pro}
\pf As earlier in this section,
let $\frakh_{\hbar}=\frakh \flb\hbar \frb$ be equipped  with the Lie bracket
$[X, Y]_\hbar =\hbar [X, Y]$, $\forall X, Y\in 
\frakh_{\hbar}$,  and $\sigma : S(\frakh  )\flb\hbar \frb  \lon
 U \frakh_{\hbar}$ the PBW-map.  Note that $\frakh_{\hbar}$ is
isomorphic to $\frakh$ as  a Lie algebra. Hence $U \frakh_{\hbar}$
is canonically isomorphic to $U \frakh \flb\hbar \frb$, whose
elements can be considered as left invariant  (formal) differential
operators on $H$.
 To each polynomial function on $T^* H\cong \frakh^* \times H$, we
assign a  (formal) differential operator on $H$ according
 to the following rule. 
For $ f\in C^{\infty} (H )$, we assign the operator   multiplying 
by $f$; for   $ f\in \mbox{Pol}(\frakh^* ) \cong S(\frakh  )$,
we assign   the left invariant differential operator $ \Vec{\sigma (f)}$;
in general, for $f(x)g(\lambda ) $ with $f(x)\in C^{\infty} (H )$
and $g(\lambda )\in \mbox{Pol}(\frakh^* )$, we assign 
the differential operator $f(x)\Vec{\sigma (g)}$.
Then the multiplication on the algebra of   differential 
operators induces an associative
  product $*_{\hbar}$ on  $\mbox{Pol}( T^* H  )\flb\hbar \frb$,
 hence a star product on $T^* H$.  It is simple to see, from the above 
construction,   that 

\begin{enumerate}
\item for any $f, g\in C^{\infty}(\frakh^* )$,
\begin{equation}
\label{eq:20}
f(\lambda )*_{\hbar}g (\lambda )=f(\lambda ) \sstar g(\lambda );
\end{equation}
\item
for any $f(x) \in  C^{\infty}(H)$ and $g(\lambda )\in C^{\infty}(\frakh^* )$,
\begin{equation} 
\label{eq:21}
f(x)*_{\hbar}g (\lambda )=f(x)g (\lambda );
\end{equation}
\item
  for any $f (\lambda )\in C^{\infty}(\frakh^* )$
and $  g  (x)\in  C^{\infty}(H)$,
\begin{equation}
\label{eq:22}
f(\lambda )*_{\hbar}g (x)
=\sum_{k=0}^{\infty} \frac{ \hbar^k }{k!}
\frac{\partial^{k} f(\lambda )}{\partial \lambda^{i_{1}}\cdots \partial \lambda^{i_{k}}}
\vh{i_{1}}\cdots \vh{i_{k}}g(x);
\end{equation}
\item
 for any $f(x), \ g(x)\in  C^{\infty}(H)$,
\begin{equation}
\label{eq:23}
f(x)*_{\hbar}g(x)= f(x)g(x).
\end{equation}
\end{enumerate}
In other words, this  is indeed  a \prefer star product
with $F\equiv 1$. Equation (\ref{eq:19}) thus follows immediately
 from Theorem \ref{thm:formula}. \qed
{\bf Remark.} It would be interesting to compare Equation (\ref{eq:19})
with the general  construction of star products
on cotangent symplectic manifolds in \cite{BNW, BNPW}.\\\\

Equation (\ref{eq:22}) implies that the element 
$f (\lambda +\hbar h)\in C^{\infty}(\frakh^* )\ot
U\frakh \flb\hbar \frb$, being considered
as  a left invariant differential operator
on $H$, admits the following expression:
$$\Vec{f (\lambda +\hbar h)}=f(\lambda )*_{\hbar}$$
Thus we have:


\begin{cor}
For any $f, g\in  C^{\infty}(\frakh^* )$,
\begin{equation}
(f \sstar g)(\lambda +\hbar h )=f(\lambda +\hbar h ) \sstar g(\lambda +\hbar h),
\end{equation}
where   the $\sstar $  on the left hand side stands for
the PBW-star product on $\frakh^* $, while
on the right hand  side it refers to the  multiplication
on the algebra  tensor product of  $(C^{\infty}(\frakh^* )\flb\hbar \frb,
 * )$ with  $U\frakh \flb\hbar \frb$.
\end{cor}
\pf Let $*_{\hbar}$ denote the star product on $T^* H$ as in Proposition
\ref{pro:H}.
For any $\phi (x)\in C^{\infty} (H)$,
\be
&&(f (\lambda )*_{\hbar }g (\lambda ))*_{\hbar } \phi (x) \ \ \mbox{(by 
Equations (\ref{eq:20}, \ref{eq:22}))}\\
&=&\sum_{k=0}^{\infty} \frac{ \hbar^k }{k!}
\frac{\partial^{k} (f (\lambda )\sstar
 g (\lambda ))}{\partial \lambda^{i_{1}}\cdots \partial \lambda^{i_{k}}}
\vh{i_{1}}\cdots \vh{i_{k}}\phi (x)\\
&=&\Vec{(f \sstar g)(\lambda +\hbar h )}\phi (x).
\ee
On the other hand,
\be
&&f (\lambda )*_{\hbar }(g (\lambda )*_{\hbar } \phi (x))
 \ \ \mbox{(by Equation (\ref{eq:19}))}\\
&=&\sum_{k=0}^{\infty} \frac{ \hbar^k }{k!}
\frac{\partial^{k} f (\lambda )}{\partial \lambda^{i_{1}}\cdots
 \partial \lambda^{i_{k}}} \sstar
\vh{i_{1}}\cdots \vh{i_{k}}(g (\lambda )*_{\hbar }\phi (x) )\\
&=&
 \sum_{k=0}^{\infty} \sum_{l=0}^{\infty} \frac{ \hbar^k }{k!}
\frac{\partial^{k} f (\lambda )}{\partial \lambda^{i_{1}}\cdots
 \partial \lambda^{i_{k}}}\sstar \vh{i_{1}}\cdots \vh{i_{k}}
(\frac{ \hbar^l}{l!}\frac{\partial^{l}
  g(\lambda )}{\partial \lambda^{j_{1}}\cdots \partial \lambda^{j_{l}}}
\vh{j_{1}}\cdots \vh{j_{l}} \phi (x))\\
&=& \sum_{k=0}^{\infty} \sum_{l=0}^{\infty}
\frac{ \hbar^{k+l} }{k!l!}
\frac{\partial^{k} f (\lambda )}{\partial \lambda^{i_{1}}\cdots
 \partial \lambda^{i_{k}}}
\sstar \frac{\partial^{l}
  g(\lambda )}{\partial \lambda^{j_{1}}\cdots \partial \lambda^{j_{l}}}
\vh{i_{1}}\cdots \vh{i_{k}}\vh{j_{1}}\cdots \vh{j_{l}} \phi (x)\\
&=&\Vec{f(\lambda +\hbar h ) \sstar g(\lambda +\hbar h )}\phi (x).
\ee
The conclusion thus follows from the associativity 
of $*_{\hbar}$. \qed

\begin{cor}
\label{cor:FG}
For any  $F, G\in C^{\infty}(\frakh^* )\ot U\frakg \ot U\frakg \flb\hbar \frb$,
\begin{equation}
(F\sstar G)_{23}( \lambda +\hbar h^{(1)})
=F_{23}( \lambda +\hbar h^{(1)})\sstar G_{23}( \lambda +\hbar h^{(1)}).
\end{equation}
In particular, if $F(\lambda )\in C^{\infty}(\frakh^* )\ot U\frakg
 \ot U\frakg \flb\hbar \frb$
is invertible,  we have
\begin{equation}
F^{-1}_{23}(\lambda +\hbar h^{(1)}) =F_{23}( \lambda +\hbar h^{(1)})^{-1}.
\end{equation}
\end{cor}

\section{Quantum dynamical Yang-Baxter equation}

The main  purpose of  this section is to derive the quantum dynamical 
Yang-Baxter equation over a nonabelian base  $\frakh$ from the
``twisted-cocycle"  condition (\ref{eq:shifted0}). This was
standard  when   $\frakh$ is Abelian (e.g., see \cite{BBB}).
The proof was  based on  the Drinfel'd  theory of
quasi-Hopf algebras \cite{dr:quasi}.
In our situation, however,
 the  quasi-Hopf algebra approach   does not work any more. 
Nevertheless, one can carry out a proof in a way completely
parallel to the ordinary case. 

The main result of this section  is the following:

\begin{thm}
\label{thm:dybe}
Assume that  $F: \frakh^* \lon U\frakg \ot U\frakg \flb\hbar \frb$
satisfies the ``twisted-cocycle" condition (\ref{eq:shifted0}).
Then  
\begin{equation}
R(\lambda )=F_{21}(\lambda )^{-1}\sstar F_{12}(\lambda )
\end{equation}
 satisfies the following generalized quantum
dynamical  Yang-Baxter equation (or  Gervais-Neveu-Felder equation):
\begin{equation}
\label{eq:dybe}
R_{12}(\lambda ) \sstar R_{13}(\lambda +\hbar h^{(2)} ) 
 \sstar R_{23}(\lambda )
=R_{23}(\lambda +\hbar h^{(1)} )  \sstar
R_{13}(\lambda )  \sstar R_{12}(\lambda +\hbar h^{(3)} ).
\end{equation}
Here  $\sstar$ denotes the natural multiplication on
$C^{\infty}(\frakh^* )\ot ( U\frakg )^{n}\flb\hbar \frb$, $\forall n$,
with $C^{\infty}(\frakh^* )$ being equipped with the PBW-star
product.
\end{thm}


It is simple to see that the  usual relation
\begin{equation}
\label{eq:ab}
\td (a \sstar b)=\td a  \sstar \td b 
\end{equation}
still holds for any $a, \ b \in C^{\infty}(\frakh^* )\ot
 U\frakg \flb\hbar \frb$. 
Define $\ttd: C^{\infty}(\frakh^* )\ot U\frakg \flb\hbar \frb\lon
C^{\infty}(\frakh^* ) \ot U\frakg \ot U\frakg \flb\hbar \frb$ by
\begin{equation}
\label{eq:30}
\ttd a =  F(\lambda )^{-1} \sstar  \td a \sstar  F(\lambda ), \ \ \ \forall
a \in C^{\infty}(\frakh^* )\ot U\frakg \flb\hbar \frb.
\end{equation}
It is  simple to see, using the associativity of $*$, that
\begin{equation}
\label{eq:31}
\ttd^{op} a =R(\lambda )\sstar \ttd a \sstar R(\lambda )^{-1}.
\end{equation}

The following is immediate  from Corollary \ref{cor:FG}.

\begin{cor}
\label{cor:R}
\begin{equation}
\label{eq:29}
R_{23}( \lambda +\hbar h^{(1)})=F_{32}(\lambda +\hbar h^{(1)})^{-1}\sstar
F_{23}( \lambda +\hbar h^{(1)}).
\end{equation}
\end{cor}
{\bf Remark.} \ \  Equation (\ref{eq:29}) is trivial when
$\frakh$ is Abelian. It, however, does  not seem obvious
in general. We can see from the proof of Corollary \ref{cor:FG}
that this equation   essentially follows from the associativity
of the   star product given by  Equation (\ref{eq:19}).\\\\

For any given $F(\lambda )\in C^{\infty}(\frakh^* )\ot U\frakg \ot  U\frakg
 \flb\hbar \frb$, introduce 
$\Phi_{123}(\lambda )\in  C^{\infty}(\frakh^* )\ot U\frakg \ot  U\frakg
\ot U\frakg \flb\hbar \frb$ by

\begin{equation}
\label{eq:33}
\Phi_{123}(\lambda )= F_{23}(\lambda )^{-1} \sstar  [(id\ot \td )F(\lambda )^{-1}]\sstar [( \td \ot  id )F(\lambda )]\sstar F_{12}(\lambda ).
\end{equation}

\begin{lem}
\begin{eqnarray}
(\ttd \ot id)R&=&\Phi_{231}\sstar R_{13}\sstar \Phi_{132}^{-1} \sstar  R_{23}
\sstar \Phi_{123}; \label{eq:34}\\
(id  \ot \ttd )R&=& \Phi_{312}^{-1} \sstar R_{13} \sstar  \Phi_{213} \sstar
R_{12} \sstar \Phi_{123}^{-1}. \label{eq:35}
\end{eqnarray}
\end{lem}
 \pf By applying the permutation $ a_{1}\ot a_{2}\ot a_{3}
\lon a_{1}\ot a_{3} \ot a_{2}$ on Equation (\ref{eq:33}), one 
obtains that 
\begin{eqnarray*}
\Phi_{132}(\lambda )&=&F_{32}(\lambda )^{-1} \sstar
  \sigma_{23} [(id\ot \td )F(\lambda )^{-1} ]\sstar 
\sigma_{23} [( \td \ot  id )F(\lambda )]\sstar F_{13}(\lambda )\\
&=&F_{32}(\lambda )^{-1} \sstar
   [(id\ot \td )F(\lambda )^{-1} ]\sstar 
\sigma_{23} [( \td \ot  id )F(\lambda )]\sstar F_{13}(\lambda ),
\end{eqnarray*}
since $\Delta$ is cocommutative.
Similarly, applying the permutation $a_{1}\ot a_{2}\ot a_{3}
 \lon  a_{2}\ot a_{3} \ot a_{1}$ on Equation (\ref{eq:33}), one
obtains that
\begin{equation}
\Phi_{231}(\lambda )=
F_{12}(\lambda )^{-1} \sstar
   [( \td  \ot id )F_{21}(\lambda )^{-1} ]\sstar
\sigma_{23} [( \td \ot  id )F(\lambda )]\sstar F_{31}(\lambda ).
\end{equation}

On the other hand,  by definition, 
\begin{eqnarray}
R_{13}(\lambda )&=&F_{31}(\lambda )^{-1} \sstar  F_{13}(\lambda )\\
R_{23}(\lambda )&=&F_{32}(\lambda )^{-1} \sstar  F_{23}(\lambda ).
\end{eqnarray}

It thus follows that 
\be
&&\Phi_{231}\sstar R_{13}\sstar \Phi_{132}^{-1} \sstar  R_{23}
\sstar \Phi_{123}\\
&=&F_{12}(\lambda )^{-1} \sstar ( \td \ot  id )F_{21}(\lambda )^{-1}
 \sstar ( \td \ot  id )F(\lambda ) \sstar  F_{12}(\lambda ) \ \ \ \mbox{(by
Equation (\ref{eq:ab}))}\\
&=&F_{12}(\lambda )^{-1} \sstar ( \td \ot  id )R(\lambda )\sstar  F_{12}(\lambda )
\ \ \ \mbox{(by Equation (\ref{eq:30}))}\\
&=& (\ttd \ot id)R. 
\ee
Equation (\ref{eq:34}) can be proved similarly. \qed
{\bf Proof of Theorem \ref{thm:dybe}.}  From Equation (\ref{eq:31}),
it follows   that 
$$R_{12}\sstar (\ttd \ot id)R =(\ttd^{op} \ot id)R \sstar R_{12}.$$
According to  Equation  (\ref{eq:34}),  this is equivalent to
$$R_{12}\sstar  \Phi_{231}\sstar R_{13}\sstar \Phi_{132}^{-1} \sstar  R_{23}
\sstar \Phi_{123}
=\Phi_{321}\sstar 
R_{23}\sstar \Phi_{312}^{-1} \sstar  R_{13}
\sstar \Phi_{213} \sstar R_{12} .$$
Thus,
\begin{equation}
R_{12}\sstar (\Phi_{231}\sstar R_{13}\sstar \Phi_{132}^{-1}) \sstar  R_{23}
=(\Phi_{321}\sstar 
R_{23}\sstar \Phi_{312}^{-1}) \sstar  R_{13}
\sstar (\Phi_{213} \sstar R_{12} \sstar \Phi_{123}^{-1} ).
\end{equation}
Now the twisted-cocycle condition (\ref{eq:shifted0}) implies that
\begin{equation}
\Phi_{123}(\lambda )=F_{12}(\lambda +
\hbar h^{(3)} )^{-1}\sstar F_{12} (\lambda ).
\end{equation}
It thus follows that
\be
&&\Phi_{213} \sstar R_{12} \sstar \Phi_{123}^{-1}\\
&=&F_{21}(\lambda + \hbar h^{(3)} )^{-1}\sstar F_{21} (\lambda )
\sstar  F_{21} (\lambda )^{-1} \sstar
F_{12} (\lambda )\sstar  F_{12} (\lambda )^{-1} \sstar
F_{12} (\lambda + \hbar h^{(3)} )\\
&=&F_{21}(\lambda + \hbar h^{(3)} )^{-1}
\sstar F_{12} (\lambda + \hbar h^{(3)}) \ \ \
 \mbox{(by Corollary \ref{cor:R}})\\
&=&R_{12} (\lambda + \hbar h^{(3)}).
\ee

Applying  the permutations: 
$a_{1}\ot a_{2}\ot a_{3} \lon  a_{3}\ot a_{1} \ot a_{2}$,
and $a_{1}\ot a_{2}\ot a_{3} \lon  a_{1}\ot a_{3} \ot a_{2}$
respectively to the equation above,
one obtains
\be
&&\Phi_{321}\sstar R_{23}\sstar \Phi_{312}^{-1} = R_{23}(\lambda + \hbar h^{(1)}) \ \ \mbox{and }\\
&&\Phi_{231}\sstar R_{13}\sstar \Phi_{132}^{-1}
=R_{13}(\lambda + \hbar h^{(2)} ).
\ee
Equation (\ref{eq:dybe}) thus follows immediately. \qed


\section{Concluding remarks}

Even though our discussion so far has been mainly
confined to  triangular dynamical r-matrices,
we should point out that there do  exist
many interesting examples of non-triangular ones.
For instance, when the Lie algebra $\frakg$ admits an ad-invariant
bilinear form and the base Lie algebra $\frakh$ equals  $\frakg$,
Alekseev and Meinrenken found an explicit construction of
an interesting non-triangular dynamical r-matrix
\cite{AM} in connection with  their study of the non-commutative Weil algebra.
In fact, for simple Lie algebras, the existence of  AM-dynamical r-matrices
was already proved by Etingof and Varchenko \cite{EV1}.
The  construction of Alekseev and Meinrenken was later generalized
by Etingof and Schiffmann to a more general context \cite{ES}.
So there is no doubt that there are abundant non-trivial
examples of  dynamical r-matrices with  a  nonabelian
base.   It is therefore    desirable
to know how they can be quantized. Inspired by the above
 discussion  in  the triangular case, we are ready to propose the
following quantization problem along the line of Drinfel'$\!$d's
naive\footnote{Drinfeld's original naive quantization
was proposed for a classical $r$-matrix in $A\ot A$ for
an associative algebra $A$. Here one can consider 
 $A$ as the universal enveloping algebra $U\frakg$,
and $r\in \frakg\ot  \frakg\subset U\frakg\ot U\frakg. $}
 quantization \cite{Drinfeld:90}.

\begin{defi}
\label{def:quantization}
Given a classical dynamical r-matrix  
  $r: \frakh^{*} \lon  \frakg \ot \frakg $,
a quantization of $r$ is $R(\lambda )=1+\hbar r(\lambda )+O(\hbar^{2} )
\in U(\frakg )\ot U(\frakg )\flb\hbar \frb$ which is $H$-equivariant and
 satisfies the  generalized quantum
dynamical  Yang-Baxter equation (or  Gervais-Neveu-Felder equation):
\begin{equation}
\label{eq:d}
R_{12}(\lambda ) \sstar R_{13}(\lambda +\hbar h^{(2)} ) 
 \sstar R_{23}(\lambda )
=R_{23}(\lambda +\hbar h^{(1)} )  \sstar
R_{13}(\lambda )  \sstar R_{12}(\lambda +\hbar h^{(3)} ).
\end{equation}
\end{defi}

Combining Proposition \ref{pro:dy}, Proposition \ref{pro:prefer},
 Corollary \ref{cor:shift} and Theorem \ref{thm:dybe},
 we may summarize the main  result  of this paper in the following:

\begin{thm}
\label{thm:main}
 A triangular dynamical $r$-matrix $r: \frakh^* \lon
\wedge^2 \frakg  $ is quantizable   if there exists a \prefer
star product  on   the corresponding Poisson
manifold  $\frakh^* \times G$.
\end{thm}

We conclude this paper with a list of questions together
 with some thoughts.

{\bf Question 1}: Is every classical  triangular 
 dynamical $r$-matrix quantizable? 

According to   Theorem \ref{thm:main}, this question
 is equivalent to
asking whether a \prefer
star product always exists for the corresponding Poisson
manifold $\frakh^* \times G$.
When the base  Lie algebra is Abelian, 
a quantization  procedure was found  for splittable classical  triangular 
 dynamical r-matrices using Fedosov's method \cite{Xu4}. 
Recently Etingof and  Nikshych \cite{EN}, using  the vertex-IRF 
transformation method, showed  the existence of  quantization  for the so called
completely degenerate triangular dynamical r-matrices,
 which leads to the hope that 
the existence of  quantization could be possibly  settled
by combing both methods in \cite{Xu4} and \cite{EN}.
However, when the
base Lie algebra $\frakh$ is nonabelian, the method in
\cite{Xu4} does not admit a straightforward
generalization. One of the  main difficulties is that the
Fedosov method uses Weyl quantization, while
our quantization here is in normal ordering. 
Nevertheless, for the dynamical  r-matrices constructed
in Theorem \ref{thm:r},  under some  mild assumptions
a quantization seems  feasible    
by using the generalized Karabegov's method 
 \cite{Ast, BD}. This  problem
 will be discussed  in a separate publication.

{\bf Question 2}: What is the symmetrized version of the quantum
dynamical Yang-Baxter equation (\ref{eq:d})?

We derived Equation (\ref{eq:d}) from a \prefer star product, which
is a normal ordering star product.  The reason for us to
choose the  normal ordering here is that one can obtain a
very explicit formula for the star product: Equation (\ref{eq:explicit}).
A Weyl ordering \prefer star product may exist, but it
may be more difficult to work with. For the canonical
cotangent  symplectic structure $T^*H$,  a 
Weyl ordering star product was  found by Gutt  \cite{Gutt},
 but it is rather difficulty to write down an   explicit  formula \cite{BB}.
 As we can  see from the previous discussion,
how  a quantum dynamical Yang-Baxter equation looks 
 is  closely related
to the   choice of  a star product on $T^* H$. When
$H$ is Abelian,  there
is a very simple operator    establishing an  isomorphism
between these two quantizations, which  is indeed the 
transformation needed to transform  a unsymmetrized QDYBE  into  a
symmetrized one.  Such an operator also exists  for a  general cotangent
bundle $T^{*}Q$ \cite{BNW}, but it is  much more  complicated.
 Nevertheless, this viewpoint  may still provide   a useful method
to obtain the symmetrized version of a QDYBE. 

{\bf Question 3}: Is every classical  dynamical $r$-matrix quantizable?

This question may be a bit too general. As a first step, it
should be already   quite interesting  to find a quantum analogue of 
Alekseev-Meinrenken  dynamical r-matrices.

{\bf Question 4}: What is the deformation theory
controlling the quantization problem as proposed  in  Definition
\ref{def:quantization}?

If  $R=1+\hbar r +\cdots +\hbar^{i} r_{i}+\cdots $,
where $r_{i} \in C^{\infty}(\frakh^* )\ot U\frakg \ot U\frakg , \ i\geq 2$,
is a solution to the QDYBE,
 the $\hbar$-term $r$ must be a solution of the
classical dynamical Yang-Baxter equation. Indeed  the
quantum dynamical Yang-Baxter
 equation implies a sequence of equations of $r_{i}$ in terms of
lower order terms. One should expect   
some cohomology theory  here just as for any deformation theory \cite{BEFFL}.
However, in our case, the equation seems very complicated.
On the other hand, it is quite surprising that such a theory 
does not  seem to exist in the literature  even  in the case of quantization
 of a   usual $r$-matrix.

Finally, we would like to point out that perhaps
a more useful way of thinking of quantization of
a dynamical r-matrix is to consider the quantum
groupoids as defined in \cite{Xu3}. This is in some
sense  an analogue of the ``sophisticated" quantization
in terms of Drinfel'd \cite{Drinfeld:90}. A classical
dynamical r-matrix gives rise to  a Lie bialgebroid
$(T\frakh^* \times \frakg, \ T^{*}\frakh^* \times \frakg^{*} )$ 
\cite{BK-S, LX1}. 
Its induced Poisson structure on the base space $\frakh^*$
is the Lie-Poisson structure $\pi_{\frakh^*}$, which admits the PBW-star 
product as a standard deformation quantization. This leads to the following

{\bf Question 5:} Does  the  Lie bialgebroid
$(T\frakh^* \times \frakg, \ T^{*}\frakh^* \times \frakg^{*} )$
corresponding to a classical dynamical r-matrix always  admit a quantization
in the sense of \cite{Xu3}, with the   base algebra being the
PBW-star algebra $C^{\infty}(\frakh^* )\flb\hbar \frb$?

To connect the quantization problem  in Definition
\ref{def:quantization} with 
that of Lie bialgebroids, it is clear that one needs to consider  preferred
quantization of Lie bialgebroids: namely, a quantization
where the total  algebra is undeformed and   remains to
be  $\cald (\frakh^* ) \ot U\frakg \flb\hbar \frb$.

{\bf Question 6:} Does  the  Lie bialgebroid
$(T\frakh^* \times \frakg, \ T^{*}\frakh^* \times \frakg^{*} )$
admit a preferred quantization?
 How is such a preferred quantization  related to
the quantization of $r$ as proposed in Definition \ref{def:quantization}?

When $\frakh=0$, namely for usual r-matrices, the answer to 
Question 6  is positive,  due
to  a remarkable theorem of Etingof-Kazhdan \cite{EK}.

\end{document}